\documentclass[11pt]{amsart}
\usepackage{mathtools}

\usepackage{pdfpages}
\usepackage[left=20mm,top=0.6in,bottom=12mm]{geometry}
\usepackage{enumitem}

\usepackage{xcolor}
\usepackage{hyperref}
\usepackage[hyperpageref]{backref}
\hypersetup{
    colorlinks=true,
    linkcolor=myblue,
    filecolor=magenta,      
    urlcolor=mygreen,
    citecolor=mygreen,
}
\definecolor{mygreen}{rgb}{0.01,0.6,0.2}
\definecolor{myblue}{rgb}{0.01, 0.18, 1.0}

\allowdisplaybreaks

\numberwithin{equation}{section}

\newtheorem{Lemma}{LEMMA}[section]
\newtheorem{Theorem}[Lemma]{Theorem}
\newtheorem{Proposition}[Lemma]{Proposition}
\newtheorem{Corollary}[Lemma]{Corollary}

\newtheorem{remark}[Lemma]{Remark}
\newtheorem{definition}[Lemma]{Definition}
\newtheorem{example}[Lemma]{Example}

\newtheorem{Fact}[Lemma]{Fact}
\newtheorem{assumption}[Lemma]{Assumption}

\def\bt{\begin{Theorem}}
\def\et{\end{Theorem}}
\def\bl{\begin{Lemma}}
\def\el{\end{Lemma}}
\def\bp{\begin{Proposition}}
\def\ep{\end{Proposition}}
\def\bcor{\begin{Corollary}}
\def\ecor{\end{Corollary}}
\def\bpf{\begin{proof}}
\def\epf{\end{proof}}

\def\brem{\begin{remark}\rm }
\def\erem{\hfill $\lozenge$ \end{remark}}

\def\bedef{\begin{definition}\rm }
\def\endef{\hfill$\lozenge$\end{definition}}

\def\beg{\begin{example}\rm }
\def\eeg{\hfill $\lozenge$\end{example}}

\def\bef{\begin{Fact}}
\def\eef{\end{Fact}}

\def\bea{\begin{assumption}}
\def\ena{\end{assumption}}
\def\bc{\begin{center}}
\def\ec{\end{center}}
\def\noi{\noindent}

\def\beq{\begin{equation}}
\def\eeq{\end{equation}}
\def\beqarray{\begin{eqnarray*}}
\def\eeqarray{\end{eqnarray*}}
\def\<{\leftangle}
\def\>{\rightangle}
\def\({\left(}
\def\){\right)}
\def\f{\varphi}

\def\<{\langle}
\def\>{\rangle}
\def\q{\quad}

\def\a{\alpha}
\def\b{\beta}

\def\d{\delta}

\def\k{\kappa}

\def\t{\tau}

\def\l{\lambda}

\def\O{\Omega}

\def\w.r.t.{with respect to}
\def\R{{\mathbb{R}}}
\def\N{{\mathbb{N}}}

\def\bq{\begin{quote}}
\def\eq{\end{quote}}

\def\bit{\begin{itemize}}
\def\eit{\end{itemize}}
\def\i{\item}
\def\ben{\begin{enumerate}}
\def\een{\end{enumerate}}

\def\ds{\displaystyle}

\begin{document}
\title[regularization for a final value problem]{A final value problem with a non-local and a source term: regularization by truncation}
\author{Subhankar Mondal}
\address{ TIFR Centre for Applicable Mathematics, Bangalore-560065, India;}
\email{subhankar22@tifrbng.res.in; 
}

\begin{abstract} 
This paper is concerned with recovering the solution of a final value problem associated with a parabolic equation involving a non linear source and a non-local term, which to the best of our knowledge has not been studied earlier. It is shown that the considered problem is ill-posed, and thus, some regularization method has to be employed in order to obtain stable approximations. In this regard, we obtain regularized approximations by solving some non linear integral equations which is derived by considering a truncated version of the Fourier expansion of the sought solution. Under different Gevrey smoothness assumptions on the exact solution, we provide parameter choice strategies and obtain the error estimates. A key tool in deriving such estimates is a version of Gr{\"o}nwalls' inequality for iterated integrals, which perhaps, is proposed and analysed for the first time.
\end{abstract}

\maketitle

\noi\textbf{Keywords:} backward problem, non-local term, ill-posed problem, regularization, parameter choice, fixed point 

\noi\textbf{MSC 2020:} 47A52, 35R25, 35R30, 45K05
\section{Introduction} 
Let $\O\subset \R^d,\,d\geq 1,$ be a bounded domain with smooth boundary $\partial\O$ and $\t>0$ be fixed. For $F:[0,\t]\times L^2(\O)\to L^2(\O)$ and $g\in L^2(\O)$ we consider the problem of recovering $u:[0,\t]\to L^2(\O)$ that satisfies the following final value problem (FVP) with a non-local term 
\beq\label{gov_pde}
\begin{cases}
u_t(x,t)-\Delta u(x,t)=F(t,u(x,t))+\int_t^\t u(x,s)\,ds,&\q \text{in}\q \O\times (0,\t),\\
u(x,t)=0,&\q \text{on}\q \partial\O\times (0,\t),\\
u(x,\t)=g(x),&\q \text{in}\q \O.
\end{cases}
\eeq
 This type of PDE appears in various context of engineering such as modelling the heat conduction for  materials with memory \cite{nunziato_1971}, analysis of space-time dependent nuclear dynamics \cite{pachpatte_1983}, epidemic phenomena in biology \cite{luan_khanh_2021}. For further details about applications of some parabolic equations with non-local term, the interested reader may refer to \cite{li_zhou_gao_2018, luan_khanh_2021} and the references therein. Additionally, it is worth to note that these type of PDE with non-local terms have attracted a lot of attention in control theory in the recent past (see for e.g., \cite{guerrero_imanuvilov_2013, ivanov_pandolfi_2009, tao_gao_2016}).

As mentioned earlier, in this paper we are interested in recovering $u:[0,\t]\to L^2(\O)$ from the knowledge of the final value $g\in L^2(\O)$ that satisfies \eqref{gov_pde} for a general source function $F:[0,\t]\times L^2(\O)\to L^2(\O).$ Throughout, we assume that the possibly non linear source function $F$ satisfies the following hypothesis:
\ben[label={$(\bf F{\arabic*})$}]
\setcounter{enumi}{0}
\i\label{F1} There exists $\k>0$ such that $\left\|F(s,\psi_1)-F(s,\psi_2)\right\|\leq \k\left\|\psi_1-\psi_2\right\|$ for all $s\in[0,\t]$ and $\psi_1,\psi_2\in L^2(\O).$ 
\i\label{F2} There exists $\phi\in L^2(\O)$ such that $s\mapsto F(s,\phi)\in L^1(0,\t;L^2(\O)).$
\een
Note that the assumptions {\rm\ref{F1}} and {\rm\ref{F2}} are quite standard. Indeed, {\rm\ref{F1}} requires that $F$ is Lipschitz \w.r.t. the second variable, and the condition {\rm\ref{F2}} is assumed just to ensure some integrability of $F$ (see Lemma \ref{prop-S-N-F} (i)).

Moreover, we assume that the considered FVP \eqref{gov_pde} has a unique solution $u\in C([0,\t];L^2(\O))$ corresponding to the exact final value $g\in L^2(\O).$ (In the next section we shall provide an example (see Example \ref{eg-sol-existence}) of a FVP of the form \eqref{gov_pde} for which we obtain explicitly the analytical representation of the solution of the corresponding FVP.)  Once we assume that a solution exists, a natural question one would like to ask is whether the solution is stable \w.r.t. the perturbation in the data. We shall show that solving the FVP \eqref{gov_pde} is an ill-posed problem (see Example \ref{eg_illposed}) in the sense that a small perturbation in the final value $g$ may lead to a large deviation in the solution $u(t)$ for each $0\leq t<\t.$ Thus, the considered FVP possesses the ill-posedness behaviour similar to that of the classical backward heat conduction equation \cite{isakov_book}.

In real world scenario, it is obvious that the final value $g$ will be available from some measurements, and hence noise in the data is inevitable. Thus, we assume that for $\d>0,$ $g^\d\in L^2(\O)$ is the noisy data corresponding to the exact final value $g$ satisfying the deterministic noise model
\beq\label{noise_level}
\left\|g^\d-g\right\|_{L^2(\O)}\leq \d.
\eeq
Thus, the problem at hand is to recover $u(t),\,0\leq t<\t,$  where $u$ is the unique solution of the FVP \eqref{gov_pde} for the exact final value $g\in L^2(\O),$ from the knowledge of $g^\d$ satisfying \eqref{noise_level}. Since the considered FVP is ill-posed, we have to employ some regularization method (see for e.g.\cite{engl_hanke_neubauer, nair_opeq}) in order to obtain some stable approximations. In this regard, we shall obtain regularized solution by solving some non linear integral equation derived using a truncated version of the Fourier expansion of the sought solution. The level at which the Fourier expansion is truncated plays the role of regularization parameter. 

The literature on the regularization aspect for the classical backward heat conduction problem is quite rich and several regularization methods (see for example, quasi reversibility method \cite{lattes_lions_1969, showalter_1974}; quasi boundary value method \cite{melnikova_1991, hao_duc_sahli_2008}; truncated regularization method \cite{elden_berntsson_reginska_2000, tuan_trong_2010, nam_2010}) were studied extensively over the last few decades. In comparison to that, there seems to be a very limited work devoted to the study of regularization aspect for the ill-posed FVP with non-local terms. To the best of our knowledge, a recent work of Luan and Khanh \cite{luan_khanh_2021} is perhaps the first work, where the authors investigated the FVP/backward problem of identification of $v(\cdot,t)$ for $0\leq t<\t$ from the knowledge of the final value $h\in L^2(\O),$ where $v$ is the solution of 
\beq\label{gov_pde_luan_khanh}
\begin{cases}
v_t(x,t)-\Delta v(x,t)=\b\int_0^t  v(x,s)\,ds,&\q \text{in}\q \O\times (0,\t),\\
v(x,t)=0,&\q \text{on}\q \partial\O\times (0,\t),\\
v(x,\t)=h(x),&\q \text{in}\q \O,
\end{cases}
\eeq
for a fixed $\b>0.$ They have shown that for $0<t<\t,$ the backward problem is well-posed in the sense that it has a unique solution and the solution depends continuously on the final value $h$, whereas for $t=0,$ the backward problem is ill-posed (see \cite[Theorem 3.2]{luan_khanh_2021}). That is, the problem of identification of $v_0:=v(\cdot,0)$ from the knowledge of the final value $h$ is ill-posed. In order to obtain stable approximations they have employed the truncated spectral regularization method and obtained error estimates of H{\"o}lder type for different type of parameter choice strategies. 

It is clear from \eqref{gov_pde} and \eqref{gov_pde_luan_khanh} that the non-local term associated with these equations are different. Also, it is worth noting that with the type of non-local term as considered in this paper, the ill-posedness behaviour of the backward problem changes completely, as discussed earlier. Thus, in view of this, studying the problem associated with \eqref{gov_pde} is indeed an interesting problem.

We now point out the important contributions of this paper. 
\bit
\i To the best of our knowledge, the FVP \eqref{gov_pde} with a general source term that allows some non-linearity and involving a different version of the non-local term is studied for the first time from the regularization point of view.
\i We establish that the considered FVP is ill-posed for every $0\leq t<\t$ and thus obtain stable approximations by employing a regularization method using the truncated Fourier expansion.
\i Under different Gevrey smoothness assumption, we obtain the corresponding error estimates. One of the key tool in deriving such rates is a version of Gr{\"o}nwalls' inequality for iterated integrals, which perhaps is proposed and analysed for the first time in this paper.  
\eit

The rest of the paper is organised as follows: In Section \ref{sec-prelim}, we recall some well known results associated to the Dirichlet eigenvalue problem for a Laplace operator. Using the eigenvalues of the Dirichlet eigenvalue problem, we obtain an integral representation of the exact solution of the FVP \eqref{gov_pde}, which plays the key role in defining the stable approximations. We also recall the definition of {\it Gevrey spaces}, that will play the role of {\it source sets}. In Section \ref{sec-prep}, we state and prove (wherever required) various necessary results that plays the essential role in deriving error estimates and defining the stable approximations. Moreover, in this section we provide a new version of Gr{\"o}nwalls' type inequality for iterated integrals. Section \ref{sec-reg} contains the main results of this paper. In this section we shall deal with regularization and obtain various error estimates. Moreover, we shall provide a parameter choice strategy and obtain the corresponding estimates.

\section{Preliminaries}\label{sec-prelim}
It is well known that for $\l\in\R$ the Dirichlet eigenvalue problem
\beq\label{eig_pde}
\begin{cases}
-\Delta v=\l v,&\q\text{in}\q\O,\\
v=0,&\q\text{on}\q\partial\O,
\end{cases}
\eeq
admits a sequence of eigenvalues $0<\l_1\leq\l_2\leq\ldots\to\infty$ and the corresponding eigenfunctions $\f_n\in H^1_0(\O)\cap H^2(\O)$ are such that $\left\{\f_n:\,n\in\N\right\}$  forms an orthonormal basis of $L^2(\O)$ (see for e.g., \cite{evans}). It is also known that (see for e.g., \cite{courant_hilbert_book}) there exist $e_1,e_2>0$ such that 
\beq\label{eigvalue_est}
e_1 n^{2/d}\leq \l_n\leq e_2n^{2/d},\q n\in\N.
\eeq

We first give an example of a FVP that is of the form \eqref{gov_pde} and show that it has a unique solution.
\beg\label{eg-sol-existence} 
Let $d=1$ and $\O=(0,1).$ Let $\l_j$ be the the eigenvalues for the Dirichlet problem \eqref{eig_pde} and $\left\{\f_j:\,j\in\N\right\}$ be the corresponding eigenfunctions. Then it is known that $\l_j=j^2\pi^2,\,j\in\N.$ Now consider the FVP
\beq\label{eg-proto}
\begin{cases}
v_t(x,t)-v_{xx}(x,t)=v(x,t)+\int_t^\t v(x,s)\,ds,&\q\text{in}\q\O\times (0,\t),\\
v(x,t)=0,&\q\text{on}\q\partial\O\times (0,\t),\\
v(x,\t)=\f_n(x),&\q\text{in}\q\O.
\end{cases}
\eeq 
Let $$\a_j=\frac{1}{2}\left(-\left(\l_j-1\right)+\sqrt{\left(\l_j-1\right)^2-4}\right)\,\,\text{and}\,\,\b_j=\frac{1}{2}\left(-\left(\l_j-1\right)-\sqrt{\left(\l_j-1\right)^2-4}\right).$$ Then it can be verified easily that $$v(x,t)=\frac{\a_ne^{-\a_n(\t-t)}-\b_ne^{-\b_n(\t-t)}}{\a_n-\b_n}\f_n(x).$$
is the solution of \eqref{eg-proto}.
\eeg

We now give an example to show that the considered FVP \eqref{gov_pde} is ill-posed in the sense that a small perturbation in the final value $g$ may lead to a large deviation in the solution $u$ of the FVP \eqref{gov_pde}. 
\beg\label{eg_illposed}
Let $\O,\l_j,\f_j\,\a_j\,\b_j$ be as considered in Example \ref{eg-sol-existence}. Consider the FVP of recovering the solution $v(\cdot,t)$ for $0\leq t<\t$ such that
\beq\label{eg-proto_illposed}
\begin{cases}
v_t-v_{xx}=v+\int_t^\t v(s)\,ds,&\q\text{in}\q\O\times (0,\t),\\
v=0,&\q\text{on}\q\partial\O\times (0,\t),\\
v(\cdot,\t)=h^n,&\q\text{in}\q\O,
\end{cases}
\eeq
where $h^n=\frac{1}{\left|\b_n\right|}\f_n.$ It can be verified easily that for every $n\in\N,$
 $$v^n(x,t)=\frac{\a_ne^{-\a_n(\t-t)}-\b_ne^{-\b_n(\t-t)}}{\left(\a_n-\b_n\right)\left|\b_n\right|}\f_n(x)$$ is the solution of \eqref{eg-proto_illposed}.
Recall that $\l_n=n^2\pi^2,\,\,\left|\a_n\right|\leq \left|\b_n\right|,\,\,\text{for all}\,\,n\in\N$ and 
$$
\left|\b_n\right|=\frac{\left(\l_n-1\right)+\sqrt{\left(\l_n-1\right)^2-4}}{2}=\frac{n^2\pi^2-1+\sqrt{\left(n^2\pi^2-1\right)^2-4}}{2}\to \infty \,\,\text{as}\,\, n\to \infty.
$$
Now, we observe that $$\left\|h^n\right\|_{L^2(\O)}=\frac{1}{\left|\b_n\right|}\to 0\q\text{as}\,\,n\to \infty.$$ Moreover, for every $0\leq t<\t,$ we have
\beqarray
\left\|v^n(t)\right\|_{L^2(\O)}&=&\left|\frac{\a_ne^{-\a_n(\t-t)}-\b_ne^{-\b_n(\t-t)}}{\left(\a_n-\b_n\right)\left|\b_n\right|}\right|=\frac{\left|\b_n\right|e^{\left|\b_n\right|(\t-t)}-\left|\a_n\right|e^{\left|\a_n\right|(\t-t)}}{\left(\left|\b_n\right|-\left|\a_n\right|\right)\left|\b_n\right|}\\
&\geq &\frac{\left|\b_n\right|e^{\left|\b_n\right|(\t-t)}-\left|\a_n\right|e^{\left|\b_n\right|(\t-t)}}{\left(\left|\b_n\right|-\left|\a_n\right|\right)\left|\b_n\right|}=\frac{e^{\left|\b_n\right|(\t-t)}}{\left|\b_n\right|},
\eeqarray
and hence, $$\left\|v^n(t)\right\|_{L^2(\O)}\to \infty\,\,\text{as}\,\,n\to \infty,\,\,0\leq t<\t.$$
This, shows that a small perturbation in the final value may lead to a large deviation in the corresponding solution of the FVP for every $0\leq t<\t.$ Thus, the problem of solving the FVP \eqref{eg-proto_illposed} is ill-posed. 
\eeg

Before proceeding further, let us fix some notations that we shall follow in the rest of the paper.
 \bit
 \i We shall denote the inner-product on $L^2(\O)$ by $\<\cdot,\cdot\>$ and the norm by $\left\|\cdot\right\|.$
 \i For a Banach space $\mathcal{X},$ $\mathcal{B}(\mathcal{X})$ shall denote the Banach space of all bounded linear operators from $\mathcal{X}$ to $\mathcal{X}.$ The norm on $\mathcal{B}\left(\mathcal{X}\right)$ will be denoted by $\left\|\cdot\right\|_{\mathcal{X}\to \mathcal{X}}.$
 \i For a Banach space $\mathcal{X},$ $C([0,\t];\mathcal{X})$ denotes the Banach space of all continuous functions $v:[0,\t]\to \mathcal{X}$ with the norm $\ds\left\|v\right\|_{C([0,\t];\mathcal{X})}:=\max_{0\leq t\leq \t} \left\|v(t)\right\|_\mathcal{X}.$
 \i $L^1(0,\t;L^2(\O))$ denotes the Banach space of all measurable functions $v:[0,\t]\to L^2(\O)$ such that $\int_0^\t \left\|v(t)\right\|\,dt<\infty$ and the norm is given by $\left\|v\right\|_{L^1(0,\t;L^2(\O))}:=\int_0^\t \left\|v(t)\right\|\,dt.$
 \i $L^\infty (0,\t;L^2(\O))$ denotes the Banach space of all measurable functions $v:[0,\t]\to L^2(\O)$ such that there exists $M>0$ satisfying $\left\|v(t)\right\|\leq M$ for almost all $t\in[0,\t],$ and the norm is given by $\left\|v\right\|_{L^\infty (0,\t;L^2(\O))}:=\inf\left\{ M:\, \left\|v(t)\right\|\leq M\,\,\text{for almost all}\,\,t\in[0,\t]\right\}.$ 
 \i For $w\in L^1(0,\t;L^2(\O))$ we shall use the standard notation $w(t)$ to denote the function $w(\cdot,t).$ Moreover, we may denote $\<w(t),\f_j\>$ by $w_j(t)$ and $\<F(t,w((t)),\f_j\>$ by $F_j(t,w(t)).$
 \i For $f\in L^2(\O)$ we shall denote $\<f,\f_j\>$ by $f_j$.
 \eit

Recall that we have assumed that the FVP \eqref{gov_pde} has a unique solution $u\in C([0,\t];L^2(\O))$ for the exact final value $g\in L^2(\O)$. Thus, taking the $L^2(\O)$ inner product on both sides of the governing equation of \eqref{gov_pde} with $\f_j$, we have for all $j\in\N,$
\beq\label{ode}
\begin{cases}
\frac{d}{dt}\left(u_j(t)\right)+\l_ju_j(t)=F_j(t,u(t))+\int_t^\t u_j(s)\,ds,&\q t\in(0,\t)\\
u_j(\t)=g_j.
\end{cases}
\eeq
Now, it can be verified easily that 
\beq\label{u_j_rep}
u_j(t)=e^{\l_j(\t-t)}g_j-\int_t^\t e^{\l_j(s-t)}F_j(s,u(s))\,ds-\int_t^\t e^{\l_j(s-t)}\int_s^\t u_j(\xi)\,d\xi\,ds.
\eeq
Thus, the solution $u$ to the FVP \eqref{gov_pde} is 
$$u(t)=\sum_{j=1}^\infty u_j(t)\f_j$$
where $u_j$ is as given in \eqref{u_j_rep}.

Moreover, for $\psi\in L^2(\O)$ if we define $\mathcal{S}(t)\psi:=\sum_{j=1}^\infty e^{\l_jt}\<\psi,\f_j\>\f_j$, then we have 
\beq\label{semigroup_rep}
u(t)=\mathcal{S}(\t-t)g-\int_t^\t \mathcal{S}(s-t)F(s,u(s))\,ds-\int_t^\t \mathcal{S}(s-t)\int_s^\t u(\xi)\,d\xi\,ds.
\eeq
The above integral equation form of the exact solution will be helpful in our analysis in the upcoming sections. In fact, we obtain the approximations for the exact solution $u$ by replacing $\mathcal{S}$ with some finite rank operators from $L^2(\O)$ to $L^2(\O).$

We end this section by recalling the definition of {\bf Gevrey space} (see for e.g., \cite{cao_rammaha_titi_1999, tuan_nguyen_au_lesnic_2017}) which shall play the role of source sets in deriving error estimates.  

For $p,q\geq 0,$ the Gevrey space $\mathbb{G}_{p,q}$ is defined as $$\mathbb{G}_{p,q}:=\left\{\psi\in L^2(\O):\sum_{j=1}^\infty \l_j^{2p}e^{2q\l_j}\left|\<\psi,\f_j\>\right|^2<\infty\right\}.$$
It is well known that $\mathbb{G}_{p,q}$ is a Hilbert space \w.r.t. the inner product
$$\<\chi,\psi\>_{\mathbb{G}_{p,q}}:=\sum_{j=1}^\infty \l_j^{2p}e^{2q\l_j}\<\chi,\f_j\>\<\f_j,\psi\>,\,\,\chi,\psi\in \mathbb{G}_{p,q},$$
so that the norm is given by $$\left\|\psi\right\|_{\mathbb{G}_{p,q}}^2=\sum_{j=1}^\infty \l_j^{2p}e^{2q\l_j}\left|\<\psi,\f_j\>\right|^2.$$

\section{Preparatory results}\label{sec-prep}
We begin by proving a version of Gr{\"o}nwalls' type inequality for iterated integrals that will play the crucial role in obtaining stability estimates in later part of the paper. This version, as mentioned in the introduction, is perhaps proposed for the first time. We should mention that this is motivated from a different version of Gr{\"o}nwalls' inequality for iterated integrals established in \cite[Theorem 1.4.1]{pachpatte_2006}.
\bl \label{gronwall_iterated_version}
Let $U:[0,\t]\to \R$ be a non-negative continuous function and $c_0,c_1>0$ be fixed. Suppose that $$U(t)\leq c_0+c_1\int_t^\t\left(U(s)+\int_s^\t U(\xi)\,d\xi\right)ds,\q\,0\leq t\leq \t.$$
Then $$U(t)\leq c_0e^{(1+c_1)(\t-t)},\,\,0\leq t\leq\t.$$
\el
\bpf
Let $$X(t):=c_0+c_1\int_t^\t\left(U(s)+\int_s^\t U(\xi)\,d\xi\right)ds,\,\,0\leq t\leq\t.$$
Then it follows that $U(t)\leq X(t)$ for all $t\in[0,\t]$ and $X(\t)=c_0.$ Thus,
\beqarray
X'(t)&=&-c_1U(t)-c_1\int_t^\t U(s)\,ds\geq -c_1\left(X(t)+\int_t^\t X(s)\,ds\right).
\eeqarray
Let $$V(t):=X(t)+\int_t^\t X(s)\,ds,\,\,t\in[0,\t].$$
Clearly, $V(\t)=X(\t)=c_0,\,0\leq X(t)\leq V(t)$ and $X'(t)\geq -c_1V(t).$ Thus,
\beqarray
V'(t)=X'(t)-X(t)\geq -c_1V(t)-V(t)=-(1+c_1)V(t).
\eeqarray
That is, $$V'(t)+(1+c_1)V(t)\geq 0, \,\,t\in(0,\t),\,\,V(\t)=c_0.$$
Thus, it follows that $$e^{(1+c_1)\t}V(\t)-e^{(1+c_1)t}V(t)\geq 0$$ and hence $$V(t)\leq c_0e^{(1+c_1)(\t-t)}.$$
Now using the relation $U(t)\leq X(t)\leq V(t),$ we have
$$U(t)\leq c_0e^{(1+c_1)(\t-t)}.$$
This completes the proof.
\epf
\brem
In Lemma \ref{gronwall_iterated_version} if we replace $c_0$ and $c_1$ by two non-negative continuous functions, say $\a(t)$ and $\b(t),$ respectively, then also one can infer similar conclusion in terms of integrals of $\a(t)$ and $\b(t).$ Since those generalizations are not needed for our purpose, we are not including that.
\erem

Recall that $g\in L^2(\O)$ is the exact final value and $u\in C([0,\t];L^2(\O))$ is the corresponding unique solution of the FVP \eqref{gov_pde}. We have seen in Example \ref{eg_illposed} that the problem of recovering $u(t)$ for \eqref{gov_pde} from the knowledge of final value $g$ is an ill-posed problem for all $0\leq t<\t.$ That is, a small perturbation in the final value may lead to a large deviation in the corresponding solution of the FVP. Since the accessible final value will be typically a measured data, noise in the final value is inevitable. Thus, it is reasonable to obtain some approximations for the exact solution from the knowledge of the measured final value/noisy data. Because of the ill-posedness of the problem, we have to employ a regularization method in order to obtain stable approximations. In this regard, we consider a Fourier truncation method for the regularization purpose. More specifically, we shall replace $\mathcal{S}$ in \eqref{semigroup_rep} by some finite rank operator described below, and then try to solve some non linear integral equation. The obtained solution will serve as the regularized approximations. Keeping this in mind, in this section we shall define those finite rank operators and study various important properties associated with it.

For $N\in\N$ and $t\in[0,\t],$ we define $\mathcal{S}_N(t):L^2(\O)\to L^2(\O)$ by
\beq\label{S-N_def}\mathcal{S}_N(t)\psi:=\sum_{j=1}^N e^{\l_jt}\<\psi,\f_j\>\f_j.\eeq
From the definition of $\mathcal{S}_N$ it is easy to verify the following result.
\bl\label{prop_S-N} 
For $N\in\N$ and $t\in[0,\t],$ let $\mathcal{S}_N(t)$ be as defined in \eqref{S-N_def}. Then the following holds:
\ben
\i[(i)] $\mathcal{S}_N(t)$ is a bounded linear operator with $\left\|S_N(t)\right\|_{L^2\to L^2}\leq e^{t\l_N}$ for every $t\in[0,\t].$
\i[(ii)] $t\mapsto \mathcal{S}_N(t)\in C([0,\t];\mathcal{B}(L^2(\O)))$
\een
\el
For $v\in C([0,\t];L^2(\O)),$ we define 
\beq\label{def-T}\mathfrak{T}v=\mathcal{S}_N(\t-t)g-\int_t^\t \mathcal{S}_N(s-t)F(s,v(s))\,ds-\int_t^\t\mathcal{S}_N(s-t)\int_s^\t v(\xi)\,d\xi\,ds.\eeq
\bl \label{prop-S-N-F}
Let $F:[0,\t]\times L^2(\O)\to L^2(\O)$ satsfies {\rm\ref{F1}} and {\rm\ref{F2}}. Let $v\in L^1(0,\t;L^2(\O)).$ Then the following hold:
\ben
\i[(i)]  $t\mapsto F(t,v(t))\in L^1(0,\t;L^2(\O)).$
\i[(ii)] $t\mapsto \int_t^\t\mathcal{S}_N(s-t)\int_s^\t v(\xi)\,d\xi\,ds\in C([0,\t];L^2(\O)).$
\i[(iii)] $t\mapsto \int_t^\t \mathcal{S}_N(s-t)F(s,v(s))\,ds\,\in C([0,\t];L^2(\O)).$
\een
\el
\bpf
Let $\k>0,\phi\in L^2(\O)$ be as in {\rm\ref{F1}} and {\rm\ref{F2}}, respectively.

\noi{\bf (i).} The proof follows by observing that
\beqarray
\int_0^\t\left\|F(s,v(s))\right\|\,ds&\leq &\int_0^\t\left(\left\|F(s,\phi)\right\|+\left\|F(s,\phi)-F(s,v(s))\right\|\right)ds\\
&\leq &\int_0^\t \left\|F(s,\phi)\right\|\,ds+\k\int_0^\t\left\|\phi-v(s)\right\|\,ds.
\eeqarray
{\bf (ii).} Let $\Psi(t):=\int_t^\t\mathcal{S}_N(s-t)\int_s^\t v(\xi)\,d\xi\,ds,\,0\leq t\leq \t.$ Then for each $t\in[0,\t]$ it follows that $\Psi(t)\in L^2(\O).$ Indeed, by Lemma \ref{prop_S-N} (i), we have
$$\left\|\Psi(t)\right\|\leq \int_t^\t \left\|\mathcal{S}_N(s-t)\int_s^\t v(\xi)\,d\xi\right\|\,ds\leq \int_t^\t e^{\l_N(s-t)}\int_s^\t\left\|v(\xi)\right\|\,d\xi\,ds.$$
We now show that $t\mapsto \Psi(t),\,0\leq t\leq\t,$ is continuous. 

Let $t_0\in[0,\t].$ Then for $t\leq t_0,$ we have
\beqarray
\Psi(t)-\Psi(t_0)&=&\int_t^\t \mathcal{S}_N(s-t)\int_s^\t v(\xi)\,d\xi\,ds-\int_{t_0}^\t\mathcal{S}_N(s-t_0)\int_s^\t v(\xi)\,d\xi\,ds\\
&=&\int_t^\t\left(\mathcal{S}_N(s-t)-\mathcal{S}_N(s-t_0)\right)\int_s^\t v(\xi)\,d\xi\,ds+\int_t^\t\mathcal{S}_N(s-t_0)\int_s^\t v(\xi)\,d\xi\,ds\\
&&\q\q-\int_{t_0}^\t \mathcal{S}_N(s-t_0)\int_s^\t v(\xi)\,d\xi\,ds\\
&=& \int_t^\t\left(\mathcal{S}_N(s-t)-\mathcal{S}_N(s-t_0)\right)\int_s^\t v(\xi)\,d\xi\,ds+\int_t^{t_0}\mathcal{S}_N(s-t_0)\int_s^\t v(\xi)\,d\xi\,ds.
\eeqarray
Therefore,
\beqarray
\left\|\Psi(t)-\Psi(t_0)\right\|&\leq &\int_t^\t \left\|\mathcal{S}_N(s-t)-\mathcal{S}_N(s-t_0)\right\|_{L^2\to L^2}\int_s^\t\left\|v(\xi)\right\|\,d\xi\,ds\\
&&\q\q+\int_t^{t_0} \left\|\mathcal{S}_N(s-t_0)\right\|_{L^2\to L^2}\int_s^\t\left\|v(\xi)\right\|\,d\xi\,ds\\
&\leq & \left(\int_0^\t\left\|\mathcal{S}_N(s-t)-\mathcal{S}_N(s-t_0)\right\|_{L^2\to L^2}\,ds+\int_t^{t_0}\left\|\mathcal{S}_N(s-t_0)\right\|_{L^2\to L^2}\,ds\right)\left\|v\right\|_{L^1(0,\t;L^2(\O))}.
\eeqarray
For $t_0\leq t,$ similarly we obtain
\beqarray
\left\|\Psi(t)-\Psi(t_0)\right\|&\leq &\left(\int_0^\t\left\|\mathcal{S}_N(s-t)-\mathcal{S}_N(s-t_0)\right\|_{L^2\to L^2}\,ds+\int_t^{t_0}\left\|\mathcal{S}_N(s-t_0)\right\|_{L^2\to L^2}\,ds\right)\left\|v\right\|_{L^1(0,\t;L^2(\O))}.
\eeqarray
Now, from Lemma \ref{prop_S-N} (ii) it follows that $\left\|\mathcal{S}_N(s-t)-\mathcal{S}_N(s-t_0)\right\|\to 0$ as $t\to t_0,$ therefore by dominated convergence theorem it follows that $\int_0^\t\left\|\mathcal{S}_N(s-t)-\mathcal{S}_N(s-t_0)\right\|\,ds\to 0$ as $t\to t_0.$ Also, $\left|\int_t^{t_0}\left\|\mathcal{S}_N(s-t_0)\right\|\,ds\right|\to 0$ as $t\to t_0.$ Thus, it follows that for each $t_0\in[0,\t],$ $\left\|\Psi(t)-\Psi(t_0)\right\|\to 0$ as $t\to t_0.$

\noi{\bf (iii).} Let $\Phi(t):=\int_t^\t \mathcal{S}_N(s-t)F(s,v(s))\,ds,\,\,0\leq t\leq \t.$ By Lemma \ref{prop_S-N} (i), it follows that $$\left\|\Phi(t)\right\|\leq \int_t^\t e^{\l_N(s-t)}\left\|F(s,v(s))\right\|\,ds.$$ Thus, by Lemma \ref{prop-S-N-F} (i) it follows that $\Phi(t)\in L^2(\O)$ for all $t\in[0,\t].$ We now show that $t\mapsto \Phi(t)$ is continuous. Let $t_0\in[0,\t].$ Then similar to the proof of Lemma \ref{prop-S-N-F} (ii), we obtain
\beqarray
\left\|\Phi(t)-\Phi(t_0)\right\|&\leq & \int_0^\t \left\|\mathcal{S}_N(s-t)-\mathcal{S}_N(s-t_0)\right\|_{L^2\to L^2}\left\|F(s,v(s))\right\|\,ds\\
&&\q\q+\left|\int_t^{t_0}\left\|\mathcal{S}_N(s-t_0)\right\|_{L^2\to L^2}\left\|F(s,v(s))\right\|\,ds\right|.
\eeqarray
 Thus, we conclude that $\left\|\Phi(t)-\Phi(t_0)\right\|\to 0$ as $t\to t_0.$
\epf

\bt\label{T_welldefined}
Let $F:[0,\t]\times L^2(\O)\to L^2(\O)$ satisfies {\rm\ref{F1}} and {\rm\ref{F2}}, and $\mathfrak{T}$ be as in \eqref{def-T}. If $v\in C([0,\t];L^2(\O)),$ then $\mathfrak{T}v\in C([0,\t];L^2(\O)).$
\et
\bpf
By Lemma \ref{prop_S-N} (ii) we observe that $t\mapsto \mathcal{S}_N(\t-t)g\in C([0,\t];L^2(\O)).$ Now the proof follows by Lemma \ref{prop-S-N-F}.
\epf
\bt\label{fixed-point}
Let $F$ and $\mathfrak{T}:C([0,\t];L^2(\O))\to C([0,\t];L^2(\O))$ be as in Theorem \ref{T_welldefined}. Then there exists a unique $v_0\in C([0,\t];L^2(\O))$ such that $\mathfrak{T}v_0=v_0.$ 
\et
\bpf
For $0\leq t\leq \t,$ we have
\beqarray
\left\|\mathfrak{T}v_1(t)-\mathfrak{T}v_2(t)\right\|&\leq & \int_t^\t\left\|\mathcal{S}_N(s-t)\left(F(s,v_1(s))-F(s,v_2(s))\right)\right\|\,ds\\
&&\q+\int_t^\t\left\|\mathcal{S}_N(s-t)\left(\int_s^\t \left(v_1(\xi)-v_2(\xi)\right)\,d\xi\right)\right\|\,ds\\
&\leq & \int_t^\t e^{\l_N(s-t)}\left(\left\|F(s,v_1(s))-F(s,v_2(s))\right\|+\int_s^\t\left\|v_1(\xi)-v_2(\xi)\right\|\,d\xi\right)\,ds\\
&\leq &\k\int_t^\t e^{\l_N(s-t)}\left\|v_1(s)-v_2(s)\right\|\,ds+ \left\|v_1-v_2\right\|_{C([0,\t];L^2(\O))}\int_t^\t e^{\l_N(s-t)}\int_s^\t \,d\xi\,ds\\
&\leq &e^{\l_N\t}\left\|v_1-v_2\right\|_{C([0,\t];L^2(\O))}\left(\k\left(\t-t\right)+\frac{\left(\t-t\right)^2}{2}\right)\\
&\leq& e^{\l_N\t}\left\|v_1-v_2\right\|_{C([0,\t];L^2(\O))}\left(\k\left(\t-t\right)+\t\left(\t-t\right)\right)
\eeqarray
and hence 
$$\left\|\mathfrak{T}v_1(t)-\mathfrak{T}v_2(t)\right\|\leq \k_0e^{\l_N\t}\left(1+\t\right)\left(\t-t\right)\left\|v_1-v_2\right\|_{C([0,\t];L^2(\O))},$$
where $\k_0=\max\{\k,1\}.$
Using this, we obtain
\beqarray
\left\|\mathfrak{T}^2v_1(t)-\mathfrak{T}^2v_2(t)\right\|&\leq &\int_t^\t\left\|\mathcal{S}_N(s-t)\left(F(s,\mathfrak{T}v_1(s))-F(s,\mathfrak{T}v_2(s))\right)\right\|\,ds\\
&&\q+\int_t^\t\left\|\mathcal{S}_N(s-t)\left(\int_s^\t \left(\mathfrak{T}v_1(\xi)-\mathfrak{T}v_2(\xi)\right)\,d\xi\right)\right\|\,ds\\
&\leq &\int_t^\t e^{\l_N(s-t)}\left(\left\|F(s,\mathfrak{T}v_1(s))-F(s,\mathfrak{T}v_2(s))\right\|+\int_s^\t\left\|\mathfrak{T}v_1(\xi)-\mathfrak{T}v_2(\xi)\right\|\,d\xi\right)ds\\
&\leq & \k \int_t^\t e^{\l_N(s-t)}\left\|\mathfrak{T}v_1(s)-\mathfrak{T}v_2(s)\right\|\,ds+\int_t^\t e^{\l_N(s-t)}\int_s^\t\left\|\mathfrak{T}v_1(\xi)-\mathfrak{T}v_2(\xi)\right\|\,d\xi\,ds\\
&\leq & \k_0e^{\l_N\t}\left(\int_t^\t \left\|\mathfrak{T}v_1(s)-\mathfrak{T}v_2(s)\right\|\,ds+\t\int_t^\t \left\|\mathfrak{T}v_1(\xi)-\mathfrak{T}v_2(\xi)\right\|\,d\xi\right)\\
&=&\k_0e^{\l_N\t}\left(1+\t\right)\int_t^\t\left\|\mathfrak{T}v_1(s)-\mathfrak{T}v_2(s)\right\|\,ds\\
&\leq & \left(\k_0e^{\l_N\t}\left(1+\t\right)\right)^2\left\|v_1-v_2\right\|_{C([0,\t];L^2(\O))}\int_t^\t \left(\t-s\right)\,ds\\
&=& \left(\k_0e^{\l_N\t}\left(1+\t\right)\right)^2\frac{\left(\t-t\right)^2}{2!}\left\|v_1-v_2\right\|_{C([0,\t];L^2(\O))}.
\eeqarray
Thus, inductively, for all $m\in\N$ we have for $0\leq t\leq \t,$
$$\left\|\mathfrak{T}^mv_1(t)-\mathfrak{T}^mv_2(t)\right\|\leq \frac{\left(\k_0e^{\l_N\t}\left(1+\t\right)\left(\t-t\right)\right)^m}{m!}\left\|v_1-v_2\right\|_{C([0,\t];L^2(\O))},$$ 
and hence
$$\left\|\mathfrak{T}^mv_1-\mathfrak{T}^mv_2\right\|_{C([0,\t];L^2(\O))}\leq \frac{\left(\k_0e^{\l_N\t}\left(1+\t\right)\t\right)^m}{m!}\left\|v_1-v_2\right\|_{C([0,\t];L^2(\O))}.$$
Since $\frac{\left(\k_0e^{\l_N\t}\left(1+\t\right)\t\right)^m}{m!}\to 0$ as $m\to \infty,$ there exists $m_0\in\N$ such that $\mathfrak{T}^{m_0}$ is a contraction. Therefore, by Banach fixed point theorem, there exists a unique $v_0\in C([0,\t];L^2(\O))$ such that $\mathfrak{T}^{m_0}v_0=v_0.$ This shows that $\mathfrak{T}\left(\mathfrak{T}^{m_0}v_0\right)=\mathfrak{T}v_0,$ that is, $\mathfrak{T}v_0$ is a fixed point of $\mathfrak{T}^{m_0}.$ Therefore, by uniqueness, we have $\mathfrak{T}v_0=v_0.$ Clearly, this $v_0$ is the unique fixed point of $\mathfrak{T}.$
\epf
\brem\label{uN-def}
Note that the unique $v_0\in C([0,\t];L^2(\O))$ as obtained in Theorem \ref{fixed-point} depends on $N$. Thus, in order to visualize that dependence, from now onwards we shall denote that $v_0$ by $u^N.$ That is, $u^N\in C([0,\t];L^2(\O))$ is the unique element such that $\mathfrak{T}u^N=u^N.$
\erem
\section{Regularization and error analysis}\label{sec-reg}
Let $u^N\in C([0,\t];L^2(\O))$ be as described in Remark \ref{uN-def}. In this section, we shall see that $u^N$ are indeed regularized approximations. Moreover, we shall obtain estimates for the quantity $\left\|u(t)-u^N(t)\right\|$ by assuming that $u$ belongs to some Gevrey spaces. Note that,  in order to obtain some estimate of $\left\|u(t)-u^N(t)\right\|$ in terms of the regularization parameter $N$ so that the estimate tends to zero as $N\to \infty,$ it is necessary to assume some apriori condition on the solution $u$ otherwise no such estimate can be derived (see for e.g., \cite{schock_1984}). Finally, we provide a parameter choice strategy, that is a procedure to choose the regularization parameter $N$, and obtain the corresponding error estimates.
\bt\label{u-uN_est}
Let $F:[0,\t]\times L^2(\O)\to L^2(\O)$ satisfies {\rm\ref{F1}} and {\rm\ref{F2}}. Let $p,q>0,$ and $\k_0=\max\{\k,1\},$ where $\k$ is as in {\rm\ref{F1}}. For $N\in\N,$ let $u^N\in C([0,\t];L^2(\O))$ be as mentioned in Remark \ref{uN-def}. Then the following hold:
\ben
\i[(i)] If $u\in L^\infty(0,\t;\mathbb{G}_{p,\t})$ then
$$\left\|u(t)-u^N(t)\right\|\leq e^{\left(1+\k_0\right)\left(\t-t\right)}\left\|u\right\|_{L^\infty(0,\t;\mathbb{G}_{p,\t})}\l_N^{-p}e^{-\l_Nt},\,\,0\leq t\leq\t.$$
\i[(ii)] If $u\in L^\infty(0,\t;\mathbb{G}_{0,q+\t})$ then $$\left\|u(t)-u^N(t)\right\|\leq e^{\left(1+\k_0\right)\left(\t-t\right)}\left\|u\right\|_{L^\infty(0,\t;\mathbb{G}_{0,q+\t})}e^{-(q+t)\l_N},\,\,0\leq t\leq \t.$$
\een
\et
\bpf
Recall that the exact solution $u$ of \eqref{gov_pde} satisfies $$u(t)=\mathcal{S}(\t-t)g-\int_t^\t\mathcal{S}(s-t)F(s,u(s))\,ds-\int_t^\t\mathcal{S}(s-t)\int_s^\t u(\xi)\,d\xi\,ds$$ and hence 
$$\sum_{j=1}^N\<u(t),\f_j\>\f_j=\mathcal{S}_N(\t-t)g-\int_t^\t\mathcal{S}_N(s-t)F(s,u(s))\,ds-\int_t^\t\mathcal{S}_N(s-t)\int_s^\t u(\xi)\,d\xi\,ds.$$
Thus,
\beqarray
\left\|u(t)-u^N(t)\right\|&\leq &\left\|u(t)-\sum_{j=1}^N u_j(t)\f_j\right\|+\left\|\sum_{j=1}^N u_j(t)\f_j-u^N(t)\right\|\\
&\leq & \sqrt{\sum_{j=N+1}^\infty \left|u_j(t)\right|^2}+\int_t^\t\left\|\mathcal{S}_N(s-t)\left(F(s,u^N(s))-F(s,u(s))\right)\right\|\,ds\\
&&\q\q+\int_t^\t\left\|\mathcal{S}_N(s-t)\left(\int_s^\t\left(u^N(\xi)-u(\xi)\right)d\xi\right)\right\|\,ds
\eeqarray
{\bf (i).} Let $u\in L^\infty(0,\t;\mathbb{G}_{p,\t}).$ Then we have
\beqarray
\left\|u(t)-u^N(t)\right\|&\leq &\sqrt{\sum_{j=N+1}^\infty \l_j^{2p}e^{2t\l_j}\l_j^{-2p}e^{-2t\l_j}\left|u_j(t)\right|^2}+\k\int_t^\t e^{\l_N(s-t)}\left\|u^N(s)-u(s)\right\|\,ds\\
&&\q\q + \int_t^\t e^{\l_N(s-t)}\int_s^\t\left\|u^N(\xi)-u(\xi)\right\|\,d\xi\,ds\\
&\leq& \l_N^{-p}e^{-t\l_N}\left\|u(t)\right\|_{\mathbb{G}_{p,\t}}\\
&&\q\q +\k_0\left(\int_t^\t e^{\l_N(s-t)}\left(\left\|u^N(s)-u(s)\right\|+\int_s^\t\left\|u^N(\xi)-u(\xi)\right\|\,d\xi\right)\,ds\right).
\eeqarray
Thus,
\beqarray
e^{t\l_N}\left\|u(t)-u^N(t)\right\|& \leq &\l_N^{-p}\left\|u\right\|_{L^\infty(0,\t;\mathbb{G}_{p,\t})}\\
&&\q+\k_0\left(\int_t^\t e^{s\l_N}\left(\left\|u^N(s)-u(s)\right\|+\int_s^\t\left\|u^N(\xi)-u(\xi)\right\|\,d\xi\right)\,ds\right).
\eeqarray
Now, note that 
\beqarray
\int_t^\t e^{s\l_N}\int_s^\t\left\|u^N(\xi)-u(\xi)\right\|\,d\xi\,ds&=&\int_t^\t e^{s\l_N}\int_s^\t e^{\xi\l_N}e^{-\xi\l_N}\left\|u^N(\xi)-u(\xi)\right\|\,d\xi\,ds\\
&\leq &\int_t^\t e^{s\l_N}\int_s^\t e^{-s\l_N}e^{\xi\l_N}\left\|u^N(\xi)-u(\xi)\right\|\,d\xi\,ds\\
&=&\int_t^\t\int_s^\t e^{\xi\l_N}\left\|u^N(\xi)-u(\xi)\right\|\,d\xi\,ds.
\eeqarray
Thus, we obtain
\beqarray
e^{t\l_N}\left\|u(t)-u^N(t)\right\|&\leq & \l_N^{-p}\left\|u\right\|_{L^\infty(0,\t;\mathbb{G}_{p,\t})}\\
&&\q\q+\k_0\left(\int_t^\t e^{s\l_N}\left\|u(s)-u^N(s)\right\|\,ds+\int_t^\t\int_s^\t e^{\xi\l_N}\left\|u(\xi)-u^N(\xi)\right\|\,d\xi\,ds\right).
\eeqarray
Therefore, by Lemma \ref{gronwall_iterated_version}, we obtain
$$e^{t\l_N}\left\|u(t)-u^N(t)\right\|\leq \l_N^{-p}\left\|u\right\|_{L^\infty(0,\t;\mathbb{G}_{p,\t})}e^{\left(1+\k_0\right)\left(\t-t\right)}.$$ 
Thus,
$$\left\|u(t)-u^N(t)\right\|\leq e^{\left(1+\k_0\right)\left(\t-t\right)}\left\|u\right\|_{L^\infty(0,\t;\mathbb{G}_{p,\t})}\l_N^{-p}e^{-\l_Nt}.$$ This completes the proof of (i).

\noi{\bf (ii).} Let $u\in L^\infty(0,\t;\mathbb{G}_{0,q+\t}).$ Then similar to the above argument we obtain
\beqarray
\left\|u(t)-u^N(t)\right\|&\leq &\sqrt{\sum_{j=N+1}^\infty e^{2(q+t)\l_j}e^{-2(q+t)\l_j}\left|u_j(t)\right|^2}+\k\int_t^\t e^{\l_N(s-t)}\left\|u^N(s)-u(s)\right\|\,ds\\
&&\q\q + \int_t^\t e^{\l_N(s-t)}\int_s^\t\left\|u^N(\xi)-u(\xi)\right\|\,d\xi\,ds\\
&\leq& e^{-(q+t)\l_N}\left\|u(t)\right\|_{\mathbb{G}_{0,q+\t}}\\
&&\q\q +\k_0\left(\int_t^\t e^{\l_N(s-t)}\left(\left\|u^N(s)-u(s)\right\|+\int_s^\t\left\|u^N(\xi)-u(\xi)\right\|\,d\xi\right)\,ds\right),
\eeqarray
and hence
\beqarray
e^{t\l_N}\left\|u(t)-u^N(t)\right\|&\leq & e^{-q\l_N}\left\|u\right\|_{L^\infty(0,\t;\mathbb{G}_{0,q+\t})}\\
&&\q\q+\k_0\left(\int_t^\t e^{s\l_N}\left\|u(s)-u^N(s)\right\|\,ds+\int_t^\t\int_s^\t e^{\xi\l_N}\left\|u(\xi)-u^N(\xi)\right\|\,d\xi\,ds\right).
\eeqarray
Therefore, by Lemma \ref{gronwall_iterated_version}, we have
$$\left\|u(t)-u^N(t)\right\|\leq e^{\left(1+\k_0\right)\left(\t-t\right)}\left\|u\right\|_{L^\infty(0,\t;\mathbb{G}_{0,q+\t})}e^{-(q+t)\l_N}.$$
\epf

We now analyse the case when we have a noisy final value in place of the exact final value. For $\d>0,$ let $g^\d\in L^2(\O)$ be the noisy final value satisfying \eqref{noise_level}.

For $v\in C([0,\t];L^2(\O))$, let 
\beq\label{def-T_delta}
\mathfrak{T}_\d v=\mathcal{S}_N(\t-t)g^\d-\int_t^\t \mathcal{S}_N(s-t)F(s,v(s))\,ds-\int_t^\t\mathcal{S}_N(s-t)\int_s^\t v(\xi)\,d\xi\,ds.
\eeq 
Then following the arguments of Theorems \ref{T_welldefined} and \ref{fixed-point} it is obvious that $\mathfrak{T}_\d:C([0,\t];L^2(\O))\to C([0,\t];L^2(\O))$ is well-defined and there exists a unique $u_\d^N\in C([0,\t];L^2(\O))$ such that 
\beq\label{uNdelta}\mathfrak{T}_\d u^N_\d=u^N_\d,\eeq that is,
$$u^N_\d(t)=\mathcal{S}_N(\t-t)g^\d-\int_t^\t\mathcal{S}_N(s-t)F(s,u^N_\d(s))\,ds-\int_t^\t\mathcal{S}_N(s-t)\int_s^\t u^N_\d(\xi)\,d\xi\,ds.$$
Having obtained $u^N_\d,$ we now obtain estimate for $\left\|u^N(t)-u^{N}_\d(t)\right\|.$
\bt\label{uN-uNdelta_est}
Let $F:[0,\t]\times L^2(\O)\to L^2(\O)$ satisfies {\rm\ref{F1}} and {\rm\ref{F2}}. Let $\k$ be as in {\rm\ref{F1}} and $\k_0=\max\{\k,1\}.$ Let $u^N, u^N_\d$ be as mentioned in Remark \ref{uN-def} and \eqref{uNdelta}, respectively. Then
$$\left\|u^N(t)-u^N_\d(t)\right\|\leq \d e^{\l_N(\t-t)}e^{(1+\k_0)(\t-t)},\q0\leq t\leq \t.$$
\et
\bpf
We have
\beqarray
\left\|u^N(t)-u^N_\d(t)\right\|&\leq &\left\|\mathcal{S}_N(\t-t)\left(g-g^\d\right)\right\|+\int_t^\t\left\|\mathcal{S}_N(s-t)\left(F(s,u^N(s))-F(s,u^N_\d(s))\right)\right\|\,ds\\
&&\q\q + \int_t^\t\left\|\mathcal{S}_N(s-t)\int_s^\t\left(u^N(\xi)-u^N_\d(\xi)\right)d\xi\right\|\,ds\\
&\leq & e^{\l_N\left(\t-t\right)}\d+\int_t^\t \k e^{\l_N(s-t)}\left\|u^N(s)-u^N_\d(s)\right\|\,ds + \int_t^\t e^{\l_N(s-t)}\int_s^\t\left\|u^N(\xi)-u^N_\d(\xi)\right\|\,d\xi\,ds\\
&\leq & e^{\l_N\left(\t-t\right)}\d\\
&&\q\q +\k_0\left(\int_t^\t e^{\l_N(s-t)}\left\|u^N(s)-u^N_\d(s)\right\|\,ds +  \int_t^\t e^{\l_N(s-t)}\int_s^\t\left\|u^N(\xi)-u^N_\d(\xi)\right\|\,d\xi\,ds\right).
\eeqarray
Thus,
\beqarray
e^{\l_Nt}\left\|u^N(t)-u^N_\d(t)\right\|&\leq & e^{\l_N\t}\d +\k_0\left(\int_t^\t e^{\l_Ns}\left\|u^N(s)-u^N_\d(s)\right\|\,ds +  \int_t^\t e^{\l_Ns}\int_s^\t\left\|u^N(\xi)-u^N_\d(\xi)\right\|\,d\xi\,ds\right)\\
&\leq & e^{\l_N\t}\d +\k_0 \left(\int_t^\t e^{\l_Ns}\left\|u^N(s)-u^N_\d(s)\right\|\,ds +  \int_t^\t \int_s^\t e^{\l_N\xi}\left\|u^N(\xi)-u^N_\d(\xi)\right\|\,d\xi\,ds\right).
\eeqarray
Therefore, by Lemma \ref{gronwall_iterated_version}, we obtain
$$e^{\l_Nt}\left\|u^N(t)-u^N_\d(t)\right\|\leq \d e^{\l_N\t}e^{(1+\k_0)(\t-t)},$$ that is,
$$\left\|u^N(t)-u^N_\d(t)\right\|\leq \d e^{\l_N(\t-t)}e^{(1+\k_0)(\t-t)}.$$
\epf
Therefore combining Theorems \ref{u-uN_est} and \ref{uN-uNdelta_est}, we have the following result.
\bt\label{final _err_est}
Let $p,q>0,$ and $\k_0=\max\{\k,1\},$ where $\k$ is as in {\rm\ref{F1}}. For $N\in\N,$ let $u^N,u^N_\d$ be as in Theorem \ref{uN-uNdelta_est}. Then the following hold:
\ben
\i[(i)] If $u\in L^\infty(0,\t;\mathbb{G}_{p,\t})$ then
$$\left\|u(t)-u^N_\d(t)\right\|\leq e^{\left(1+\k_0\right)\left(\t-t\right)}\left\|u\right\|_{L^\infty(0,\t;\mathbb{G}_{p,\t})}\l_N^{-p}e^{-\l_Nt} + \d e^{\l_N(\t-t)}e^{(1+\k_0)(\t-t)},\q0\leq t\leq \t.$$
\i[(ii)] If $u\in L^\infty(0,\t;\mathbb{G}_{0,q+\t})$ then
$$\left\|u(t)-u^N_\d(t)\right\|\leq e^{\left(1+\k_0\right)\left(\t-t\right)}\left\|u\right\|_{L^\infty(0,\t;\mathbb{G}_{0,q+\t})}e^{-(q+t)\l_N} + \d e^{\l_N(\t-t)}e^{(1+\k_0)(\t-t)},\q0\leq t\leq \t.$$
\een
\et
We now provide parameter choice strategies, that is, a strategy to choose the regularization parameter $N$ depending on $\d$ in such a way that $N\to \infty$ and $\left\|u(t)-u^N_\d(t)\right\|\to 0$ as $\d\to 0.$ We first recall an important result from \cite{tautenhahn_1998} that will be useful in the upcoming analysis.
\bl\label{tautenhahn}\cite[Lemma 3.3]{tautenhahn_1998}
Let $0<a<1, b,c,d$ be positive constants. Let $\zeta:(0,a]\to\R$ be given by $\zeta(s):=s^b\left(d\ln\frac{1}{s}\right)^{-c}.$ Then for the inverse function $\zeta^{-1}(s),$ we have
$$\zeta^{-1}(s)=s^{1/b}\left(\frac{d}{b}\ln\frac{1}{s}\right)^{c/b}\left(1+o(1)\right)\q\text{for}\,\,s\to 0.$$
\el
For $p,q,\varrho>0$, let 
\beq\label{K_p_tau}
\mathcal{K}_{p,\t}=\left\{v\in C([0,\t];L^2(\O)):\,\left\|v\right\|_{L^\infty(0,\t;\mathbb{G}_{p,\t})}\leq \varrho\right\}
\eeq and 
\beq\label{M_0_q}
\mathcal{M}_{0,q}=\left\{v\in C([0,\t];L^2(\O)):\,\left\|v\right\|_{L^\infty(0,\t;\mathbb{G}_{0,q+\t})}\leq \varrho\right\}.
\eeq
 Also, let  $\k_1=1+\k_0,$ where $\k_0$ is as in Theorem \ref{final _err_est}. Let $u\in\mathcal{K}_{p,\t}.$ Let $e_1,e_2>0$ be as in \eqref{eigvalue_est}. Then from Theorem \ref{final _err_est}, we obtain
\beqarray
\left\|u(t)-u^N_\d(t)\right\|&\leq &e^{\k_1(\t-t)}\left(\frac{\varrho}{\l_N^pe^{\l_Nt}}+\d e^{\l_N(\t-t)}\right)\\
&\leq & E_1e^{\k_1(\t-t)}\left(\frac{\varrho}{N^{2p/d}\left(e^{N^{2/d}}\right)^{e_1 t}}+\d \left(e^{N^{2/d}}\right)^{e_2(\t-t)}\right),
\eeqarray
where $E_1=\max\{e_1^{-p},1\}.$ For $r>0,$ let $A(r)=\frac{\varrho}{r^{2p/d}}+\d \left(e^{r^{2/d}}\right)^{\eta(t)},$ where $\eta(t)=e_1t+e_2(\t-t),\,\,0\leq t\leq \t.$ Note that $A(r)$ will attain minimum at $r_\d$ if $\frac{\varrho}{{r_\d}^{2p/d}}=\d\left(e^{{r_\d}^{2/d}}\right)^{\eta(t)}$ that is, if $\frac{\d}{\varrho}=\left(\frac{1}{r_\d}\right)^{2p/d}e^{-\eta(t){r_\d}^{2/d}}.$ Let $s_\d=e^{-{r_\d}^{2/d}}.$ Then it follows that $$\frac{\d}{\varrho}={s_\d}^{\eta(t)}\left(\ln\frac{1}{s_\d}\right)^{-p}.$$
Thus, by Lemma \ref{tautenhahn}, we have
$$s_\d=\left(\frac{\d}{\varrho}\right)^{1/\eta(t)}\left(\frac{1}{\eta(t)}\ln\frac{\varrho}{\d}\right)^{p/\eta(t)}\left(1+o(1)\right)\q \text{for}\,\,\d\to 0.$$
Thus,
$${r_\d}^{2/d}=\ln\left(\left(\frac{\d}{\varrho}\right)^{-1/\eta(t)}\left(\frac{1}{\eta(t)}\ln\frac{\varrho}{\d}\right)^{-p/\eta(t)}\right)\left(1+o(1)\right)\q\text{for}\,\,\d\to 0$$
and hence
$$r_\d=\left(\ln\left(\left(\frac{\d}{\varrho}\right)^{-1/\eta(t)}\left(\frac{1}{\eta(t)}\ln\frac{\varrho}{\d}\right)^{-p/\eta(t)}\right)\right)^{d/2}\left(1+o(1)\right)\q\text{for}\,\,\d\to 0.$$ 
For $r\in\R$, let $\left[\left [r\right]\right]$ denotes the greatest integer not exceeding $r.$ We define $N_\d=\left[\left[r_\d\right]\right].$ Then for this choice of the regularization parameter, from Theorem \ref{final _err_est}, it follows that for each $t\in[0,\t]$ there exists $\d_t>0$ such that 
\beqarray
\left\|u(t)-u^{N_\d}_\d(t)\right\|&\leq & 2E_1e^{\k_1(\t-t)}\,\d\left(\frac{\d}{\varrho}\right)^{-e_2\left(\t-t\right)/\eta(t)}\left(\frac{1}{\eta(t)}\ln\frac{\varrho}{\d}\right)^{-pe_2\left(\t-t\right)/\eta(t)},\q\text{for}\,\d\leq \d_t\\
&=& 2 E_1 e^{\k_1(\t-t)}\,\varrho^{e_2\left(\t-t\right)/\eta(t)}\,\d^{e_1 t/\eta(t)}\left(\frac{1}{\eta(t)}\ln\frac{\varrho}{\d}\right)^{-pe_2\left(\t-t\right)/\eta(t)},\q\text{for}\,\d\leq \d_t
\eeqarray
Now suppose that $u\in \mathcal{M}_{0,q}.$ Let $e_1,e_2>0$ be as in \eqref{eigvalue_est}. Then from Theorem \ref{final _err_est} (ii), we have
\beqarray
\left\|u(t)-u^N_\d(t)\right\|&\leq & e^{\k_1(\t-t)}\left(\frac{\varrho}{e^{(q+t)\l_N}}+\d e^{\l_N(\t-t)}\right)\\
&\leq & e^{\k_1(\t-t)}\left(\frac{\varrho}{e^{(q+t)e_1N^{2/d}}}+\d e^{(\t-t)e_2N^{2/d}}\right).
\eeqarray
For $r>0,$ let $A(r)=\frac{\varrho}{e^{e_1(q+t)r^{2/d}}}+\d e^{e_2(\t-t)r^{2/d}}.$ Clearly $A(r)$ attains it minimum at $r_\d$ if $$\frac{\varrho}{\d}=\left(e^{{r_\d}^{2/d}}\right)^{e_1(q+t)+e_2(\t-t)},$$ that is, if $$r_\d=\left(\frac{1}{e_1(q+t)+e_2(\t-t)}\ln\frac{\varrho}{\d}\right)^{d/2}.$$ Let $N_\d=\left[\left[r_\d\right]\right].$ Then from Theorem \ref{final _err_est} (ii), it follows that for each $t\in[0,\t]$ 
\beqarray
\left\|u(t)-u^{N_\d}_\d(t)\right\|&\leq & 2 e^{\k_1(\t-t)}\,\d\left(\frac{\varrho}{\d}\right)^{e_2(\t-t)/\left(e_1(q+t)+e_2(\t-t)\right)}\\
&=&2 e^{\k_1(\t-t)}\,\varrho^\frac{e_2(\t-t)}{e_1(q+t)+e_2(\t-t)}\,\d^\frac{e_1(q+t)}{e_1(q+t)+e_2(\t-t)}.
\eeqarray
Thus, we have proved the following result.
\bt
For $p,q,\varrho>0,$ let $\mathcal{K}_{p,\t}$ and $\mathcal{M}_{0,q}$ be as defined in \eqref{K_p_tau} and \eqref{M_0_q}, respectively. Let $e_1,e_2$ be as in \eqref{eigvalue_est}. Let $\k_0=\max\{\k,1\},\,\,\k_1=1+\k_0,$ where $\k$ is as in {\rm\ref{F1}}, and $\eta(t)=e_1 t+e_2(\t-t),\,\,0\leq t\leq \t.$ For $r\in\R,$ let $\left[\left[r\right]\right]$ denotes the greatest integer not exceeding $r.$ Then the following holds:
\ben
\i[(i)] If $u\in\mathcal{K}_{p,\t}$ then for each $t\in[0,\t]$ there exists $\d_t>0$ such that for the choice $$N_\d=\left[\left[\left(\ln\left(\left(\frac{\d}{\varrho}\right)^{-1/\eta(t)}\left(\frac{1}{\eta(t)}\ln\frac{\varrho}{\d}\right)^{-p/\eta(t)}\right)\right)^{d/2}\right]\right],$$ we have
$$\left\|u(t)-u^{N_\d}_\d(t)\right\|=O\left(e^{\k_1(\t-t)}\varrho^{e_2(\t-t)/\eta(t)}\,\d^{e_1t/\eta(t)}\left(\frac{1}{\eta(t)}\ln\frac{\varrho}{\d}\right)^{-pe_2\left(\t-t\right)/\eta(t)}\right),\q \d\to 0.$$
\i[(ii)] If $u\in \mathcal{M}_{0,q}$ then for each $t\in[0,\t]$ and $$N_\d=\left[\left[\left(\frac{1}{e_1(q+t)+e_2(\t-t)}\ln\frac{\varrho}{\d}\right)^{d/2}\right]\right],$$ we have
$$\left\|u(t)-u^{N_\d}_\d(t)\right\|=O\left(e^{\k_1(\t-t)}\,\varrho^\frac{e_2(\t-t)}{e_1(q+t)+e_2(\t-t)}\,\d^\frac{e_1(q+t)}{e_1(q+t)+e_2(\t-t)}\right),\q\d\to 0.$$
\een
\et
\brem
From the above Theorem it follows that if $u\in\mathcal{K}_{p,\t}$ then for $t=0$ we obtain logarithmic rate of convergence and for $0<t< \t$ we obtain a logarithmic-type of convergence. Moreover, if $u\in\mathcal{M}_{0,q}$ then we obtain a H{\"o}lder rate of convergence for every $t\in[0,\t].$
\erem
\section{Conclusion}
We have considered the problem of recovering the solution of a final value problem for a parabolic equation with a non-linear source and a non-local term, from the knowledge of the final value. We have shown by an example that the considered problem is ill-posed, that is, a small perturbation in the data may lead to a large deviation in the sought solution. Since the problem is ill-posed, we have obtained regularized solutions by solving a non-linear integral equation which is derived by considering the truncated version of Fourier representation of the solution. Under appropriate Gevrey smoothness assumption (i.e., source condition), we have established a logarithmic-type convergence rate for the regularized approximation and with some comparatively higher Gevrey smoothness assumption, H{\"o}lder rate of convergence was established for the regularized approximations. We have proposed and proved a new version of Gr{\"o}nwalls' inequality for iterated integrals, which plays the key role in obtaining the above mentioned rates.

\vspace{1cm}
\noi {\bf Acknowledgements.} The author is supported by the post doctoral fellowship of TIFR Centre for Applicable Mathematics, Bangalore.

\end{document}